\documentclass[reqno,11pt]{amsart}

\newtheorem{theorem}{Theorem}[section]

\newtheorem{lemma}{Lemma}[section]

\theoremstyle{definition}

\newtheorem{remark}{Remark}

\numberwithin{equation}{section}
\renewcommand{\b}[1]{\mbox{\boldmath $#1$}}
\newcommand{\lef}{\langle\hskip -1.8pt \langle}
\newcommand{\rig}{\rangle\hskip -1.8pt \rangle}

\DeclareMathOperator*{\trace}{trace}

\begin{document}
% ----------------------------------------------------------------

\author{A. Guillin}
\address{CEREMADE, Universit\'e Paris Dauphine}
\email{guillin@ceremade.dauphine.fr}

\author{R.~Liptser}
\address{Department of Electrical Engineering-Systems,
Tel Aviv University, 69978 Tel Aviv Israel}
\email{liptser@eng.tau.ac.il}

\title[]
{\bf MDP for integral functionals of fast and slow processes with
averaging}

\subjclass{60J27,60F10}
\keywords{Moderate deviations, Poisson equation, Puhalskii theorem}
\thanks{}

\begin{abstract}
We establish large deviation principle (LDP) for the family of
vector-valued random processes
$(X^\varepsilon,Y^\varepsilon),\varepsilon\to 0$ defined as
$$
\begin{aligned}
& X^\varepsilon_t=\frac{1}{\varepsilon^\kappa}\int_0^t
H(\xi^\varepsilon_s,Y^\varepsilon_s)ds,
\\
& dY^\varepsilon_t=F(\xi^\varepsilon_t,Y^\varepsilon_t)dt+
              \varepsilon^{1/2-\kappa}G(\xi^\varepsilon_t,Y^\varepsilon_t)dW_t,
\end{aligned}
$$
where $W_t$ is Wiener process and $\xi^\varepsilon_t$ is fast
ergodic diffusion. We show that, under $\kappa<\frac{1}{2}$ or
less  and Veretennikov-Khasminskii type condition for fast diffusion, the LDP holds
with rate function of Freidlin-Wentzell's type.
\end{abstract}
\maketitle

\section{Introduction}
\label{sec-1}

In this paper, we examine a large deviation principle (LDP) for a
family of integral functionals (here, $H(z,y)$ is a vector-valued
function of the size $p$)
\begin{equation}\label{D1}
X^\varepsilon_t=
\frac{1}{\varepsilon^\kappa}\int_0^tH(\xi^\varepsilon_s,Y^\varepsilon_s)ds,
\ \varepsilon\to  0, \ 0<\kappa< 1/2\text{(or less)}
\end{equation}
of arguments ``$\xi^\varepsilon,Y^\varepsilon$'' the diffusion
pair with fast $\xi^\varepsilon_t\in\mathbb{R}^d$ and slow
$Y^\varepsilon_t\in\mathbb{R}^\ell$ entries that solve It\^o's
equations (with respect to independent standard vector Wiener processes
$B_t$ and $W_t$)
\begin{equation}\label{1.2m}
\begin{aligned}
d\xi^\varepsilon_t&=\frac{1}{\varepsilon}b(\xi^\varepsilon_t,Y^\varepsilon_t)dt+
             \frac{1}{\sqrt{\varepsilon}}\sigma(\xi^\varepsilon_t,Y^\varepsilon_t)dB_t
\\
dY^\varepsilon_t&=F(\xi^\varepsilon_t,Y^\varepsilon_t)dt+
              \varepsilon^{1/2-\kappa}G(\xi^\varepsilon_t,Y^\varepsilon_t)dW_t
\end{aligned}
\end{equation}
subject to fixed initial conditions $\xi^\varepsilon_0=z_0$,
$Y^\varepsilon_0=y_0$.

A choice of ``$\kappa<1/2$'' is imposed by two reasons. As
follows, e.g. from Pardoux and Veretennikov \cite{PVII},
$\kappa=1/2$ defines {\it central limit theorem scale} and so,
only a diffusion approximation for the family
$X^\varepsilon,Y^\varepsilon$, $\varepsilon\to 0$ might be
expected with a limit $X,Y$ where $X$ is a continuous martingale
and $Y$ a diffusion process. The ``$0\le\kappa<1/2$'' defines {\it
large deviation principle scale} in which ``$\kappa=0$'' is at
most difficult for LDP  analysis without any guaranty of existing
an explicit formula for the rate function; the latter makes this
LDP useless for applications. The proposed ``$0<\kappa<1/2$(or
less \footnote{The meaning ``or less'' is that
$\kappa<\left(1-\frac{m}{2}\right) \wedge\frac{1}{2}$, where
parameter $m\in (0,2)$  reflects dependence of slow and fast
diffusions, see assumption (A$_{\kappa+m}$).})'', keeping in {\it
large deviations principle scale}, if the LDP holds, possesses the
rate function of Freidlin-Wentzell's type. An explanation of this
phenomena follows from {\it fast} approximation
\begin{equation}\label{appx}
X^\varepsilon_t=\varepsilon^{1/2-\kappa}
\frac{1}{\sqrt{\varepsilon}}\int_0^tH(\xi^\varepsilon_s,Y^\varepsilon_s)ds
\stackrel{\text{fast}}{\approx}\varepsilon^{1/2-\kappa}
\int_0^th(\xi^\varepsilon_s,Y^\varepsilon_s)dB_s,
\end{equation}
for some $h(z,y)$, and {\it fast} homogenization for (here $^*$ is the
transposition symbol)
\begin{gather*}
\int_0^thh^*(\xi^\varepsilon_s, Y^\varepsilon_s)ds
\\
\int_0^tF(\xi^\varepsilon_s, Y^\varepsilon_s)ds, \ \
\int_0^tGG^*(\xi^\varepsilon_s, Y^\varepsilon_s)ds.
\end{gather*}

If what was told is correct, ``$0<\kappa<\left(1-\frac{m}{2}\right)\wedge 1/2$'' is said
to belong
to {\it moderate deviation principle scale}. In this paper, we
will follow this terminology and say MDP (moderate deviation
principle) instead of LDP.

Let us refer to models close to considered here. We mention first
the case when $\xi^\varepsilon_t\equiv \xi_{t/\varepsilon}$ with
$$
d\xi_t=b(\xi_t)dt+\sigma(\xi_t)dB_t,
$$
and is assumed to satisfy Veretennikov-Khasminskii's condition (see (A$_a$),
(A$_b$) in Section \ref{sec-ass}).
If, moreover,
$H(z,y)\equiv H(z)$ and $ \int_{\mathbb{R}^d}H(z)\pi(z)dz=0, $
where $\pi(z)$ is the invariant density of fast process, the MDP
for $X^\varepsilon$ is known from Wu \cite{wu1} and Guillin
\cite{G1}, \cite{G2}. For $p=1$ and $H(z,y)=H_1(z)H_2(y)$, the MDP
type evaluation (upper bound) is given in Liptser and Spokoiny
\cite{LSpok1}. With same fast process, the MDP for $Y^\varepsilon$
can be found in Liptser, Spokoiny and Veretennikov \cite{LSV}
under $F(z,y)\equiv F(z)$. A helpful role in verification of MDP
in \cite{LSpok1} and \cite{LSV} plays the Poisson equation
$$
\mathcal{L}u(z)=-H(z)
$$
with $\mathcal{L}$ the diffusion operator of $\xi$. A role of
Poisson equation in an examination of functional central limit
theorem based on a method of corrector is well known from
Papanicolaou, Stroock and Varadhan \cite{PSV}, Ethier and Kurtz
\cite{EK}, Bhattacharya \cite{Br}, Pardoux and Veretennikov
\cite{PVI}, \cite{PVII}. In the examination of MDP this equation
plays the same role allowing to check \eqref{appx}. In this paper,
the Poisson equation implementation is borrowed from Pardoux and
Veretennikov papers \cite{PVI}, \cite{PVII}.

Thus, we deal with MDP for family of continuous random processes
even having one smooth component. The MDP evaluation results are
well known for many different settings: Borovkov, Mogulski
\cite{BM1}, \cite{BM2} and
 Chen \cite{Ch}, Ledoux \cite{Led} (processes with independent increments);
Dembo \cite{D} (martingales with bounded jumps); Dembo and Zajic
\cite{DZa} (functional of empirical processes); Dembo and Zeitouni
\cite{DZst} (iterates of expanding maps); Puhalskii \cite{Pque}
(queues in critical loading); Chang, Yao, Zajic \cite{Zaj} (queues
with long-range dependent input); Wu \cite{wu}, de Acosta, Chen
\cite{dAC}, Guillin \cite{G0}, Djellout, Guillin \cite{DG} (Markov
chains).

In the setting of present paper, following \cite{PVII}, we
introduce a diffusion process $z^y$
parametrized by $y\in\mathbb{R}^\ell$:
\begin{equation}\label{diffu}
dz^y_t=b(z^y_t,y) dt + \sigma(z^y_t,y)dB_t,
\end{equation}
subject to $z^y_0=z_0$. The ergodic property of $z^y_t$ is
provided by Veretennikov-Khasminskii type condition
(A$_a$)+(A$_b$) (Section \ref{sec-ass}, for more details see
Khasminskii \cite{Kh} and Veretennikov, \cite{Ver98}). With
$\mathcal{L}^y$ the diffusion operator $z^y$, we use the Poisson
equation
\begin{equation}\label{eqP}
\mathcal{L}^yu(z,y) = - H(z,y)
\end{equation}
for a verification of \eqref{appx} and also of validity for fast
homogenization effect.

The paper is organized as follow. In Section \ref{sec-ass},
assumptions are given and the MDP is formulated. In Section
\ref{sec-prel}, \eqref{appx} is substantiated. The proof of MDP is
given in Section \ref{sec-4}. All auxiliary results are gathered
in Appendix.

\section{Notations, assumptions and main result}
\label{sec-ass}

Hereafter, the following notations are used.

Denote $\lef\cdot,\cdot\rig$ and  $\|\cdot\|$ the inner product and
Euclidean norm respectively.

Set
$$
a(z,y):=\sigma\sigma^*(z,y)\quad\text{and}\quad A(z,y):=GG^*(z,y).
$$
The fact that $f(z,y)$ is partially differentiable (with bounded
partial derivatives) $i$-times in $z$, $j$-times in $y$ is
indicated by $\partial^i_z\partial^j_yf(z,y)$; also
$\partial^i_z\partial^j_yf(z,y)$ denotes the partial derivative
itself.

A generic nonnegative increasing function growing not faster than
polynomial is denoted by $g=g(v), v\in\mathbb{R}_+$.

\bigskip
The first group of assumptions is concerned to functions
$a,b,F,G$:

\medskip
$\b{(A_a)}$ $\partial^3  _z\partial^2_ya(z,y)$; for some
$\Lambda>\lambda>0$ and $I$ the identical matrix
$$
\lambda I\leq a(x,y)\leq \Lambda I.
$$

\medskip
$\b{(A_b)}$ $\partial^3  _z\partial^2_yb(z,y)$; for some $r>0$ and $C>0$
$$
\sup_{y}\lef z,
b(z,y)\rig\le -r\|z\|^2, \ \|z\|\ge C.
$$

\medskip
$\b{(A_F)}$ $\partial^2_yF(z,y)$; $F$ is Lipschitz continuous and
$\|F(z,y)\|\le K(1+\|y\|)$.

\medskip
$\b{(A_G)}$ $\partial^2_yG(z,y)$; $G$ is Lipschitz continuous and
bounded.

\begin{remark}
From Pardoux and Veretennikov \cite{PVI}, \cite{PVII}, it follows
that, under $\b{(A_a)}$ and $\b{(A_b)}$, the diffusion process
$z^y=(z^y_t)_{t\ge 0}$, defined in \eqref{diffu}, is ergodic with
the unique invariant density $\pi(z;y)$.
\end{remark}

\bigskip
The second group of assumptions is concerned to function $H$.

\medskip
$\b{(A_H)}$ $\partial^2_yH(z,y)$;

\hskip 10pt 1. $\|H(z,y)\|\le
g(\|y\|)(1+\|z\|)^{-(2+\varepsilon)}$, for  some small
$\varepsilon>0$

\hskip 10pt 2. $\|\partial_y H(z,y)\|+\|\partial^2_y H(z,y)\|\le
g(\|y\|)$

\hskip 10pt 3. $\int_{\mathbb{R}^d}H(z,y)\pi(z;y)dz\equiv 0$.

\bigskip
The third group of assumptions is concerned with solution of
Poisson equation \eqref{eqP} with
\begin{equation}\label{Ly}
\mathcal{L}^y=
\frac{1}{2}\sum_{i,j=1}^{d}(\sigma\sigma^*)_{ij}(z,y)\frac{\partial^2}
{\partial z_i\partial z_j}
+\sum_{i=1}^{d}b_i(z,y)\frac{\partial}{\partial z_i}.
\end{equation}

\begin{remark}
It is known from Pardoux and Veretennikov, \cite{PVI} (Proposition
1, Theorem 2), \cite{PVII} (Theorems 3) that $u(z,y) $ the
centered solution of \eqref{eqP},
$$
\int_{\mathbb{R}^d}u(z,y)\pi(z;y)dz\equiv 0,
$$
possesses the following properties:

\begin{equation}\label{bound-dery}
\begin{aligned}
&\|u(z,y)\|+\|\nabla_z u(z,y)\|\le g(\|y\|)
\\
&\exists \ m_1>0 \ \text{such that:}
\\
&\|\partial_y u(z,y)\|+\|\partial^2_y u(z,y)\|\le
g(\|y\|)(1+\|z\|^{m_1}).
\end{aligned}
\end{equation}
Also, it can be derived from Theorem 9, \cite{PVII}, that
\begin{equation}\label{bound-dery2}
\begin{aligned}
&\exists \ m_2>0 \ \text{such that:}
\\
& \|\partial^2_{zy} u(z,y)\|+\|\partial^3_{zyy} u(z,y)\|\le
g(\|y\|)(1+\|z\|^{m_2}).
\end{aligned}
\end{equation}
\eqref{bound-dery2}.
\end{remark}

Set
\begin{equation}\label{2.5sfQ}
\begin{aligned}
& Q(z,y)=\nabla_zu(z,y)\sigma\sigma^*(z,y)\big(\nabla_zu
(z,y)\big)^*
\\
&\mathsf{Q}(y)=\int_{\mathbb{R}^d}Q(z,y)\pi(z;y)dz
\\
&\mathsf{A}(y)=\int_{\mathbb{R}^d}A(z,y)\pi(z;y)dz
\\
&\mathsf{F}(y)=\int_{\mathbb{R}^d}F(z,y)\pi(z;y)dz.
\end{aligned}
\end{equation}

\bigskip
$\b{(A_{A+Q})}$ $\mathsf{A}(y)$ and $\mathsf{Q}(y)$ are uniformly
nonsingular matrices.

\medskip
Let $\mathsf{V}(z,y)$ denotes any of functions
$$
F(z,y)-\mathsf{F}(y), \ \ A(z,y)-\mathsf{A}(y), \ \
Q(z,y)-\mathsf{Q}(y)
$$
and let $V(z,y)$ denotes any entries of $\mathsf{V}(z,y)$. For $n$
large enough, choose
\begin{equation}\label{2.7m}
V'(z,y)=
  \begin{cases}
    V(z,y), & \|z\|\le n \\
    0, & \|z\|>n+1
  \end{cases}
\quad\text{with}\quad \int_{\mathbb{R}^d}V'(z,y)\pi(z;y)dz\equiv 0
\end{equation}
and set
\begin{equation}\label{V''}
V''(z,y)=V(z,y)-V'(z,y).
\end{equation}

Consider the Poisson equation with $H$ replaced by $V'$:
\begin{equation}\label{eqPV}
\mathcal{L}^y\mathfrak{u}(z,y)=-V'(z,y).
\end{equation}
By (A$_F$), (A$_G$), (A$_a$) and \eqref{bound-dery},
\eqref{bound-dery2}, the function $V'$ obeys the properties of
$H$ mentioned in (A$_H$), so that the function $\mathfrak{u}$
possesses the properties of $u$ mentioned in \eqref{bound-dery};
particularly,
\begin{equation}\label{bound-dery3}
\begin{aligned}
&\|\mathfrak{u}(z,y)\|+\|\nabla_z \mathfrak{u}(z,y)\|\le g(\|y\|)
\\
&\exists \ m_3>0 \ \text{such that:}
\\
&\|\partial_y \mathfrak{u}(z,y)\|+\|\partial^2_y \mathfrak{u}
(z,y)\|\le g(\|y\|)(1+\|z\|^{m_3}).
\end{aligned}
\end{equation}

\bigskip
$\b{(A_{\kappa+m})}$ For $m=\max(m_1,m_2,m_3)$,

\hskip 10pt 1. $m<2$

\hskip 10pt 2. $\kappa<\Big(1-\frac{m}{2}\Big)\wedge\frac{1}{2}.$

\begin{remark}
To avoid many intricate computations and help reader to follow
easily proofs, we restrict ourselves by a little-bit less general
assumptions than those which might be sufficient.
\end{remark}

\medskip
For family $(X^\varepsilon,Y^\varepsilon)$, defined in \eqref{D1},
\eqref{1.2m}, we examine the MDP in
$\mathbb{C}_{[0,\infty)}(\mathbb{R}^p\times\mathbb{R}^\ell)$
supplied by the local uniform metric
$$
\rho_\infty\big(\{X',Y'\},\{X'',Y''\}\big)=\sum\limits_{n=1}^\infty
2^{-n} \Big(1\wedge \rho_n\big(\{X',Y'\},\{X'',Y''\}\big)\Big),
$$
where for any $T>0$
$$
\begin{aligned}
& \rho_T\big(\{X',Y'\},\{X'',Y''\}\big)
\\
&=\sup\limits_{t\le T}\Bigg(\sum\limits_{i=1}^p|X'_t(i)-X''_t(i)|+
\sum\limits_{i=1}^\ell|Y'_t(i)-Y''_t(i)|\Bigg)
\end{aligned}
$$
and $X'_t(i),X''_t(i)$ and  $Y'_t(i),Y''_t(i)$ are entries of
$X'_t,X''_t$ and $Y'_t,Y''_t$ respectively.

We follow the standard Varadhan's definition of the LDP,
\cite{V}.

Set
\begin{equation}\label{JXY}
\mathsf{J}(X,Y)=
  \begin{cases}
    \frac{1}{2}\int_0^\infty\|\dot{X}_t\|^2_{\mathsf{Q}^{-1}(Y_t)}+
    \|\dot{Y}_t-\mathsf{F}(Y_t)\|^2_{\mathsf{A}^{-1}(Y_t)}, &
    {\substack{X_0=0,Y_0=y_0\\
    dX_t=\dot{X}_tdt,dY_t=\dot{Y}_tdt}}
    \\
    \infty, & \text{otherwise}.
  \end{cases}
\end{equation}

\begin{theorem}\label{theo-mt}
Under (A$_a$), (A$_b$), (A$_F$), (A$_G$), (A$_H$), (A$_{A+Q}$),
(A$_{\kappa+m}$), the family
$(X^\varepsilon,Y^\varepsilon),\varepsilon\to 0$, obeys the MDP in
the metric space
$$
\Big(\mathbb{C}_{[0,\infty)}(\mathbb{R}^p\times\mathbb{R}^\ell), \
\rho_\infty\Big)
$$
with the rate of speed $\varepsilon^{1-2\kappa}$ and rate function
given in \eqref{JXY}.
\end{theorem}

\section{Preliminaries}
\label{sec-prel}

It this section, we present new family
$(\widehat{X}^\varepsilon,Y^\varepsilon)$, with
$
\widehat{X}^\varepsilon_t=\varepsilon^{1/2-\kappa}M^\varepsilon_t
$
and
\begin{equation}\label{mart}
M^\varepsilon_t=\int_0^t\nabla_zu
(\xi^\varepsilon_s,
Y^\varepsilon_s)\sigma(\xi^\varepsilon_s,Y^\varepsilon_s)dB_s,
\end{equation}
which shares, announced in Theorem
\ref{theo-mt}, MDP  with $(X^\varepsilon,Y^\varepsilon)$
provided that at least for one of these families this MDP holds true.
To this end, it suffices to show that $\Delta^\varepsilon_t:=X^\varepsilon_t-
\widehat{X}^\varepsilon_t$
obeys the following property: for any $T>0$ and $\eta>0$,
\begin{equation}\label{disting}
\lim_{\varepsilon\to 0}\varepsilon^{1-2\kappa}\log
P\Big(\sup_{t\le T}\|\Delta ^\varepsilon_t\|>\eta\Big)=-\infty.
\end{equation}

A key role in a verification of \eqref{disting} plays Poisson equation \eqref{eqP}.
\begin{lemma}\label{lem-3.1}
Under the assumption of Theorem \ref{theo-mt}, \eqref{disting}
holds true.
\end{lemma}
\begin{proof}
Parallel to the diffusion operator $\mathcal{L}^y$ of $(z^y_t)$ we
introduce also diffusion operator
$$
\mathfrak{L}^\varepsilon
=\frac{\varepsilon^{1-2\kappa}}{2}\sum_{i,j}
(GG^*)_{ij}(z,y)\frac{\partial^2}{\partial y_i\partial y_j}+
\sum_iF_i(z,y)\frac{\partial}{\partial y_i}.
$$
With a help of $\mathcal{L}^y$ and $\mathfrak{L}^\varepsilon$, the
It\^o formula, applied to $u(\xi^\varepsilon_t,Y^\varepsilon_t)$,
gives
\begin{equation}\label{3.2mx}
\begin{aligned}
u(\xi^\varepsilon_t,Y^\varepsilon_t)&=u(z_0,y_0)+\frac{1}{\varepsilon}
\int_0^t\mathcal{L}^yu (\xi^\varepsilon_s,Y^\varepsilon_s)ds
\\
&\quad +\frac{1}{\sqrt{\varepsilon}}\int_0^t\nabla_zu
(\xi^\varepsilon_s, Y^\varepsilon_s)
\sigma(\xi^\varepsilon_s,Y^\varepsilon_s)dB_s
+\int_0^t\mathfrak{L}^\varepsilon u
(\xi^\varepsilon_s,Y^\varepsilon_s)ds
\\
&\quad +\varepsilon^{1/2-\kappa}\int_0^t\nabla_yu
(\xi^\varepsilon_s,Y^\varepsilon_s)
G(\xi^\varepsilon_s,Y^\varepsilon_s)dW_s.
\end{aligned}
\end{equation}
By \eqref{eqP}, $ \mathcal{L}^yu
(\xi^\varepsilon_s,Y^\varepsilon_s)=-H(\xi^\varepsilon_s,Y^\varepsilon_s)
$ and so, by \eqref{3.2mx} and \eqref{mart}, we have
$$
\begin{aligned}
&\frac{1}{\varepsilon^\kappa}\int_0^tH(\xi^\varepsilon_s,Y^\varepsilon_s)ds
=\widehat{X}^\varepsilon_t
\\
&\quad +\Bigg\{\varepsilon^{1-\kappa}\big[u (z_0,y_0)-u
(\xi^\varepsilon_t,Y^\varepsilon_t)\big]
+\varepsilon^{1-\kappa}\int_0^t\mathfrak{L}^\varepsilon u
(\xi^\varepsilon_s,Y^\varepsilon_s)ds
\\
&\qquad\qquad+\varepsilon^{3/2-2\kappa}\int_0^t\nabla_yu
(\xi^\varepsilon_s,Y^\varepsilon_s)
G(\xi^\varepsilon_s,Y^\varepsilon_s)dW_s\Bigg\}.
\end{aligned}
$$
Now, owing to $X^\varepsilon_t=
\frac{1}{\varepsilon^\kappa}\int_0^tH(\xi^\varepsilon_s,Y^\varepsilon_s)ds$, we get the
following presentation for $\Delta^\varepsilon_t$:
$$
\begin{aligned}
\Delta^\varepsilon_t&=\varepsilon^{1-\kappa}\big[u(z_0,y_0)
-u(\xi^\varepsilon_t,Y^\varepsilon_t)\big]
+\varepsilon^{1-\kappa}\int_0^t\mathfrak{L}^\varepsilon
u(\xi^\varepsilon_s,Y^\varepsilon_s)ds
\\
&\quad +\varepsilon^{3/2-2\kappa}\int_0^t\nabla_y
u(\xi^\varepsilon_s,Y^\varepsilon_s)
G(\xi^\varepsilon_s,Y^\varepsilon_s)dW_s.
\end{aligned}
$$
Obviously, \eqref{disting} is valid, if for any $T>0$ and $\eta>0$
\begin{equation}\label{neg1}
\begin{aligned}
& \lim_{\varepsilon\to 0}\varepsilon^{1-2\kappa}\log
P\Big(\varepsilon^{1-\kappa}\sup_{t\le T}\big\|
u(\xi^\varepsilon_t,Y^\varepsilon_t)\big\|>\eta\Big)=-\infty
\\
& \lim_{\varepsilon\to 0}\varepsilon^{1-2\kappa}\log
P\Big(\varepsilon^{1-\kappa}\int_0^T\big\|\mathfrak{L}^\varepsilon
u(\xi^\varepsilon_s,Y^\varepsilon_s)\big\|ds>\eta\Big) =-\infty
\end{aligned}
\end{equation}
and
\begin{equation}\label{neg2}
\lim_{\varepsilon\to 0}\varepsilon^{1-2\kappa}\log
P\Big(\sup_{t\le T}\varepsilon^{3/2-2\kappa}\Big\|
\int_0^t\nabla_yu(\xi^\varepsilon_s,Y^\varepsilon_s)
G(\xi^\varepsilon_s,Y^\varepsilon_s)dW_s\Big\|>\eta\Big) =-\infty.
\end{equation}

By \eqref{bound-dery},
$$
\sup_{t\le T}\big\|
u(\xi^\varepsilon_t,Y^\varepsilon_t)\big\|\le g(\sup_{t\le T}\|Y^\varepsilon_t\|)
$$
and by $\b{(A_F)}$, $\b{(A_G)}$, \eqref{bound-dery} and $\b{(A_{\kappa+m})}$,
$$
\sup_{t\le T}\big\|\mathfrak{L}^\varepsilon
u(\xi^\varepsilon_t,Y^\varepsilon_t)\big\|\le g(\sup_{t\le T}\|Y^\varepsilon_t\|)
(1+\sup_{t\le T}\|\xi^\varepsilon_t\|^m).
$$
Hence, \eqref{neg1} holds true provided that

\begin{equation}\label{3.9m}
\lim_{\varepsilon\to 0}\varepsilon^{1-2\kappa}\log
P\Big(\varepsilon^{1-\kappa}\sup_{t\le
T}\|\xi^\varepsilon_t\|^m>\eta\Big) =-\infty,
\end{equation}
and
\begin{equation}\label{3.6mx}
\lim_{C\to\infty}\varlimsup_{\varepsilon\to
0}\varepsilon^{1-2\kappa}\log P\Big( g(\sup_{t\le
T}\|Y^\varepsilon_t\|)>C\Big) =-\infty.
\end{equation}
Notice that \eqref{3.9m} follows from Lemma \ref{theo-A.1} (Appendix \ref{App-A})
with $l=1-k$ and $p=m$ with $m$ from (A$_{\kappa+m}$).
The validity of \eqref{3.6mx}, due to
$$
\Big\{g(\sup_{t\le T}\|Y^\varepsilon_t\|)>C\Big\}=
\Big\{\sup_{t\le T}\|Y^\varepsilon_t\|>g^{-1}(C)\Big\},
$$
is checked with a help of
Lemma \ref{lem-A} (Appendix \ref{sec-7}) (the assumptions of Lemma \ref{lem-A}
are provided by $(A_F)$, $(A_G)$).
A verification of \eqref{neg2} heavily uses the fact that
$
N^\varepsilon_t=\int_0^t\nabla_yu(\xi^\varepsilon_s,Y^\varepsilon_s)
G(\xi^\varepsilon_s,Y^\varepsilon_s)dW_s
$
is a continuous martingale with the
predictable quadratic variation process
$$
\langle N^\varepsilon\rangle_t=\int_0^t\nabla_yu(\xi^\varepsilon_s,Y^\varepsilon_s)
GG^*(\xi^\varepsilon_s,Y^\varepsilon_s)
\big(\nabla_yu(\xi^\varepsilon_s,Y^\varepsilon_s)\big)^*ds.
$$
It clear that it suffices to prove \eqref{neg2} for any entry of $N^\varepsilon_t$.
Let $n^\varepsilon_t$ denote any entry of $N^\varepsilon_t$ and
$\langle n^\varepsilon\rangle_t$ be its predictable variation process.

Write
\begin{multline*}
P\Big(\sup_{t\le T}\varepsilon^{3/2-2\kappa}|
n^\varepsilon_t|>\eta\Big)
= P\Big(\sup_{t\le T}\varepsilon^{1/2-\kappa}|
\varepsilon^{1-\kappa} n^\varepsilon_t|>\eta\Big)
\\
\le 2\Bigg(P\Big(\sup_{t\le T}|\varepsilon^{1-\kappa}
n^\varepsilon_t|>\frac{\eta}{\varepsilon^{1/2-\kappa}},\varepsilon^{2-2\kappa}
\langle n^\varepsilon\rangle_T\le \delta\Big)
\\
 \bigvee P\Big(\varepsilon^{1-\kappa}\langle
n^\varepsilon\rangle^{1/2}_T> \delta^{1/2}\Big) \Bigg).
\end{multline*}
By Lemma \ref{lem-B.1}(3) (Appendix \ref{App-B}), it holds
$$
P\Big(\sup_{t\le T}|\varepsilon^{1-\kappa}
n^\varepsilon_t|>\frac{\eta}{\varepsilon^{1/2-\kappa}},\varepsilon^{2-2\kappa}
\langle n^\varepsilon\rangle_T\le \delta\Big)\le 2\exp
\Big(-\frac{\eta^2}{2\varepsilon^{(1-2\kappa)}\delta^2}\Big),
$$
while, by $(A_G)$ and \eqref{bound-dery}, $ \langle
n^\varepsilon\rangle_T \le (1+\sup_{t\le
T}\|\xi^\varepsilon_t\|^m)^2g^2(\sup_{t\le T}\|Y^\varepsilon_t\|).
$

Then
$$
\begin{aligned}
&\varlimsup_{\varepsilon\to 0}\varepsilon^{1-2\kappa} \log
P\Big(\sup_{t\le T}\varepsilon^{1/2-\kappa}|
\varepsilon^{1-\kappa} n^\varepsilon_t|>\eta\Big)
\\
&\le-\frac{\eta^2}{2\delta}\bigvee \varlimsup_{\varepsilon\to
0}\varepsilon^{1-2\kappa} \log
P\Big(\varepsilon^{1-\kappa}(1+\sup_{t\le
T}\|\xi^\varepsilon_t\|^m) g(\sup_{t\le
T}\|Y^\varepsilon_t\|)>\delta^{1/2}\Big).
\end{aligned}
$$
From \eqref{3.9m} and \eqref{3.6mx}, it is readily to derive
$$
\varlimsup_{\varepsilon\to
0}\varepsilon^{1-2\kappa} \log
P\Big(\varepsilon^{1-\kappa}(1+\sup_{t\le
T}\|\xi^\varepsilon_t\|^m) g(\sup_{t\le
T}\|Y^\varepsilon_t\|)>\delta^{1/2}\Big)=-\infty.
$$
Consequently,
$$
\varlimsup_{\varepsilon\to 0}\varepsilon^{1-2\kappa} \log
P\Big(\sup_{t\le T}\varepsilon^{1/2-\kappa}|
\varepsilon^{1-\kappa} n^\varepsilon_t|>\eta\Big)\to
-\frac{\eta^2}{2\delta}\to-\infty,  \ \delta\to0.
$$
\end{proof}

\section{The proof of Theorem \ref{theo-mt}}
\label{sec-4} By Lemma \ref{lem-3.1}, it suffices to verify the MDP, announced in
Theorem \ref{theo-mt}, for family $(\widehat{X}^\varepsilon_t,Y^\varepsilon_t)$.

Recall that
\begin{equation}\label{4.1m}
\begin{aligned}
&
Y^\varepsilon_t=y_0+\int_0^tF(\xi^\varepsilon_s,Y^\varepsilon_s)ds+
\varepsilon^{1/2-\kappa}\int_0^tG(\xi^\varepsilon_s,Y^\varepsilon_s)dW_s
\\
& \widehat{X}^\varepsilon_t=\varepsilon^{1/2-\kappa}\int_0^t\nabla_zu
(\xi^\varepsilon_s,
Y^\varepsilon_s)\sigma(\xi^\varepsilon_s,Y^\varepsilon_s)dB_s.
\end{aligned}
\end{equation}

With functions $\mathsf{F}(y)$, $\mathsf{A}(y)$, and
$\mathsf{Q}(y)$, defined in \eqref{2.5sfQ}, let us introduce
``averaged analog '' of \eqref{4.1m}:
$$
\begin{aligned}
&
y^\varepsilon_t=y_0+\int_0^t\mathsf{F}(y^\varepsilon_s)ds+
\varepsilon^{1/2-\kappa}\mathsf{A}^{1/2}(y^\varepsilon_t)dW_t
\\
&
x^\varepsilon_t=\varepsilon^{1/2-\kappa}\int_0^t\mathsf{Q}^{1/2}
(y^\varepsilon_s)dB_s.
\end{aligned}
$$
Since matrices
$\mathsf{A}$ and $\mathsf{Q}$ are uniformly nonsingular, by
Freidlin and Wentzell, \cite{FW}, the family $
(x^\varepsilon_t,y^\varepsilon_t)_{t\ge 0}, \
\varepsilon\to 0 $ possesses the LDP (in terminology of this
paper, MDP) announced in Theorem \ref{theo-mt}. So, it remains to
prove that families $(\widehat{X}^\varepsilon_t,Y^\varepsilon_t)$ and
$(x^\varepsilon_t,y^\varepsilon_t)$ share this LDP.

Since matrices $\mathsf{A}$ and $\mathsf{Q}$ are uniformly
nonsingular, Puhalskii's Theorem 2.3, \cite{puh1}, (see
also \cite{puh2}), adapted to the case considered,
 is the most convenient tool for such verification. Following this
theorem, we need to verify the fast homogenization for
$F(\xi^\varepsilon,Y^\varepsilon_t)$,
$A(\xi^\varepsilon,Y^\varepsilon_t)$ and
$Q(\xi^\varepsilon,Y^\varepsilon_t)$: for any $T>0$ and $\eta>0$
\begin{equation}\label{4.3m}
\lim_{\varepsilon\to  0}\varepsilon^{1-2\kappa}\log
P\Big(\sup_{t\le T}
\Big\|\int_0^t\mathsf{V}(\xi^\varepsilon_s,Y^\varepsilon_s)ds\Big\|>\eta\Big)=-\infty,
\end{equation}
where $\mathsf{V}(z,y)$ denotes any of $Q(z,y)-\mathsf{Q}(y)$,
$F(z,y)-\mathsf{F}(y)$, $A(z,y)-\mathsf{A}(y)$. Obviously, is
suffices to verify \eqref{4.3m} only for any entries  of
$\mathsf{V}(z,y)$, denoted in Section \ref{sec-ass} by $V(z,y)$:
$$
\lim_{\varepsilon\to  0}\varepsilon^{1-2\kappa}\log
P\Big(\sup_{t\le T}
\Big|\int_0^tV(\xi^\varepsilon_s,Y^\varepsilon_s)ds\Big|>\eta\Big)=-\infty.
$$
Recall that $V=V'+V''$, where $V'$ and $V''$ are defined in
\eqref{2.7m} and \eqref{V''} respectively.

\begin{lemma}\label{lem-4.1m}
Under the assumptions of Theorem \ref{theo-mt}, for any $T>0$,
$\eta>0$
$$
\begin{aligned}
&\rm{(i)}\quad \lim_{\varepsilon\to  0}\varepsilon^{1-2\kappa}\log
P\Big(\sup_{t\le T}
\Big|\int_0^tV'(\xi^\varepsilon_s,Y^\varepsilon_s)ds\Big|>\eta\Big)=-\infty
\\
&{\rm (ii)}\quad \lim_{\varepsilon\to
0}\varepsilon^{1-2\kappa}\log P\Big(\sup_{t\le T}
\Big|\int_0^tV''(\xi^\varepsilon_s,Y^\varepsilon_s)ds\Big|>\eta\Big)=-\infty.
\end{aligned}
$$
\end{lemma}
\begin{proof}
(i). With $\mathfrak{u}$, being solution of Poisson equation
\eqref{eqPV},
 let us define a random process
$\mathfrak{u}(\xi^\varepsilon_t,Y^\varepsilon_t)$. By the
It\^o formula
\begin{multline}\label{4.2mx}
\mathfrak{u}(\xi^\varepsilon_t,Y^\varepsilon_t)=\mathfrak{u}
(z_0,y_0)+\frac{1}{\varepsilon} \int_0^t\mathcal{L}^y\mathfrak{u}
(\xi^\varepsilon_s,Y^\varepsilon_s)ds
\\
\qquad+\frac{1}{\sqrt{\varepsilon}}\int_0^t\nabla_z\mathfrak{u}
(\xi^\varepsilon_s, Y^\varepsilon_s)
\sigma(\xi^\varepsilon_s,Y^\varepsilon_s)dB_s
+\int_0^t\mathfrak{L}^\varepsilon \mathfrak{u}
(\xi^\varepsilon_s,Y^\varepsilon_s)ds
\\
+\varepsilon^{1/2-\kappa}\int_0^t\nabla_y\mathfrak{u}
(\xi^\varepsilon_s,Y^\varepsilon_s)
G(\xi^\varepsilon_s,Y^\varepsilon_s)dW_s.
\end{multline}
From \eqref{4.2mx} and \eqref{eqPV}, it
follows
$$
\begin{aligned}
\int_0^tV'(\xi^\varepsilon_s,Y^\varepsilon_s)ds&=
\sqrt{\varepsilon}\int_0^t\nabla_z\mathfrak{u} (\xi^\varepsilon_s,
Y^\varepsilon_s) \sigma(\xi^\varepsilon_s,Y^\varepsilon_s)dB_s
\\
&\quad+\Bigg\{\varepsilon[\mathfrak{u}
(z_0,y_0)-\mathfrak{u}(\xi^\varepsilon_t,Y^\varepsilon_t)]
+\varepsilon\int_0^t\mathfrak{L}^\varepsilon \mathfrak{u}
(\xi^\varepsilon_s,Y^\varepsilon_s)ds
\\
&\qquad\qquad+\varepsilon^{3/2-\kappa}\int_0^t\nabla_y\mathfrak{u}
(\xi^\varepsilon_s,Y^\varepsilon_s)
G(\xi^\varepsilon_s,Y^\varepsilon_s)dW_s\Bigg\}\\
&:=\sqrt{\varepsilon}\widehat{M}^\varepsilon_t+\widehat{\Delta}^\varepsilon_t,
\end{aligned}
$$
where
$$
\widehat{M}^\varepsilon_t=\int_0^t\nabla_z\mathfrak{u}
(\xi^\varepsilon_s, Y^\varepsilon_s)
\sigma(\xi^\varepsilon_s,Y^\varepsilon_s)dB_s
$$
and
$$
\begin{aligned}
\widehat{\Delta}^\varepsilon_t&=\Bigg\{\varepsilon[\mathfrak{u}
(z_0,y_0)-\mathfrak{u}(\xi^\varepsilon_t,Y^\varepsilon_t)]+
\\
&\quad+\varepsilon\int_0^t\mathfrak{L}^\varepsilon \mathfrak{u}
(\xi^\varepsilon_s,Y^\varepsilon_s)ds+
\varepsilon^{3/2-\kappa}\int_0^t\nabla_y\mathfrak{u}
(\xi^\varepsilon_s,Y^\varepsilon_s)
G(\xi^\varepsilon_s,Y^\varepsilon_s)dW_s\Bigg\}.
\end{aligned}
$$
It suffices to show that
\begin{equation}\label{4.7m}
\begin{aligned}
&\rm{1.}\quad \lim_{\varepsilon\to  0}\varepsilon^{1-2\kappa}\log
P\big(\sup_{t\le T}
\big|\widehat{\Delta}^\varepsilon_t\big|>\eta\big)=-\infty
\\
&{\rm 2.}\quad \lim_{\varepsilon\to  0}\varepsilon^{1-2\kappa}\log
P\big(\sup_{t\le T}
\big|\sqrt{\varepsilon}\widehat{M}^\varepsilon_t\big|>\eta\big)=-\infty.
\end{aligned}
\end{equation}
The functions $u$ and $\mathfrak{u}$, being solutions of Poisson's equations
\eqref{eqP} and \eqref{eqPV} respectively, by
\eqref{bound-dery3} and \eqref{bound-dery} possess similar properties. Moreover, since
$\varepsilon<\varepsilon^{1-\kappa}$ and
$\varepsilon^{3/2-\kappa}<\varepsilon^{3/2-2\kappa}$, the proof
of 1. in \eqref{4.7m} is similar to  the proof of
Lemma \ref{lem-3.1}.

To check 2. in \eqref{4.7m}, notice that $\sqrt{\varepsilon}
\widehat{M}^\varepsilon_t$ is a continuous martingale with the
predictable quadratic variation process
$$
\langle
\sqrt{\varepsilon}\widehat{M}^\varepsilon\rangle_t=\varepsilon
\int_0^t\nabla_z\mathfrak{u} (\xi^\varepsilon_s,
Y^\varepsilon_s)\sigma(\xi^\varepsilon_s,Y^\varepsilon_s)
\sigma^*(\xi^\varepsilon_s,Y^\varepsilon_s) (\nabla_z\mathfrak{u}
(\xi^\varepsilon_s, Y^\varepsilon_s))^*ds.
$$
Since $\big\langle\sqrt{\varepsilon}\widehat{M}^\varepsilon\big\rangle_t=
\varepsilon\big\langle \widehat{M}^\varepsilon\big\rangle_t$,
write
$$
\begin{aligned}
P\big(\sup_{t\le
T}\big|\sqrt{\varepsilon}\widehat{M}^\varepsilon_t\big|>\eta\big)&\le
2\Big\{P\big(\sup_{t\le
T}\big|\sqrt{\varepsilon}\widehat{M}^\varepsilon_t\big|>\eta,
\big\langle \sqrt{\varepsilon} M^\varepsilon\big\rangle_T\le \varepsilon
C\big)
\\
&\bigvee P\big(\big\langle \widehat{M}^\varepsilon\big\rangle^{1/2}_T>C^{1/2}\big)\Big\}.
\end{aligned}
$$
Hence,
$$
\begin{aligned}
&
\varlimsup_{\varepsilon\to 0}\varepsilon^{1-2\kappa}\log P\big(\sup_{t\le
T}\big|\sqrt{\varepsilon}\widehat{M}^\varepsilon_t\big|>\eta\big)
\\
&\quad\le
\varlimsup_{\varepsilon\to 0}\varepsilon^{1-2\kappa}\log P\big(\sup_{t\le
T}\big|\sqrt{\varepsilon}\widehat{M}^\varepsilon_t\big|>\eta,
\big\langle \sqrt{\varepsilon} M^\varepsilon\big\rangle_T\le \varepsilon
C\big)
\\
&\hskip .5in \bigvee\varlimsup_{\varepsilon\to 0}\varepsilon^{1-2\kappa}\log
P\big(\big\langle \widehat{M}^\varepsilon\big\rangle^{1/2}_T>C^{1/2}\big)\Big\}.
\end{aligned}
$$
By Lemma \ref{lem-B.1} 3), we have
$$
\varepsilon^{1-2\kappa}\log P\big(\sup_{t\le
T}\big|\sqrt{\varepsilon}\widehat{M}^\varepsilon_t\big|>\eta,
\big\langle \sqrt{\varepsilon} M^\varepsilon\big\rangle_T\le \varepsilon
C\big)\le \varepsilon^{1-2\kappa}\log 2-\frac{\eta^2}{2C\varepsilon^{2\kappa}}.
$$
Consequently,
$$
\varlimsup_{\varepsilon\to 0}\varepsilon^{1-2\kappa}\log P\big(\sup_{t\le
T}\big|\sqrt{\varepsilon}\widehat{M}^\varepsilon_t\big|>\eta\big)
\le\varlimsup_{\varepsilon\to 0}\varepsilon^{1-2\kappa}\log
P\big(\big\langle \widehat{M}^\varepsilon\big\rangle^{1/2}_T>C^{1/2}\big).
$$
Notice now that by (A$_a$) and \eqref{bound-dery3},
$\big\langle \widehat{M}^\varepsilon\big\rangle^{1/2}_T\le g(\sup_{t\le T}
\|Y^\varepsilon_t\|)$. Hence, by \eqref{3.6mx},
$
\lim_{C\to\infty}\varlimsup_{\varepsilon\to 0}\varepsilon^{1-2\kappa}\log
P\big(\big\langle \widehat{M}^\varepsilon\big\rangle^{1/2}_T>C^{1/2}\big)\to-\infty.
$

Thus, the second part in \eqref{4.7m} is valid.

\medskip
(ii) From the definition of $V''(z,y)$ (see \eqref{V''}), it
follows that
$$
|V''(z,y)|\le g(y)I(\|z\|>n).
$$
Therefore, the desired statement holds provided that
$$
\lim_{\varepsilon\to  0}\varepsilon^{1-2\kappa}\log
P\Big(g(\sup_{t\le T} \|Y^\varepsilon_t\|)
\int_0^TI(\|\xi^\varepsilon_s\|>n)ds>\eta\Big)=-\infty.
$$
On the other hand, taking into account \eqref{3.6mx}, it suffices
to prove only that for any $T>0$ and $\eta>0$
\begin{equation}\label{I>n}
\lim_{\varepsilon\to  0}\varepsilon^{1-2\kappa}\log P\Big(
\int_0^TI(\|\xi^\varepsilon_s\|>n)ds>\eta\Big)=-\infty.
\end{equation}

A verification of \eqref{I>n} uses a nonlinear operator
$\mathcal{D}^y$, introduced in Liptser \cite{L} (see (4.16) there,
and also applied in Liptser, Spokoiny and Veretennikov \cite{LSV},
Section 4.1.2), acting on $v=v(z)$  the twice continuously
differentiable function as follows:
\begin{equation}\label{Dy}
\mathcal{D}^yv(z)=\mathcal{L}^yv(z)+\frac{1}{2}\|\nabla_z
v(z)\sigma(z,y)\|^2.
\end{equation}
We apply $\mathcal{D}^y$ to $ v(z)=\frac{\|z\|^2}{1+\|z\|}$ and
notice that the gradient of this function is defined as
$$
\nabla_z v(z)=\frac{\|z\|(2+\|z\|)}{(1+\|z\|)^2}\frac{z}{\|z\|}
$$
and is bounded. This property, jointly with (A$_a$), provides the
boundedness of $\|\nabla_z v(z)\sigma(z,y)\|^2$. Also, the
boundedness of $\partial^2v$ is readily verified.

Set
$$
U_t^\varepsilon=v(\xi^\varepsilon_t)-v(z_0)
-\frac{1}{\varepsilon}\int_0^t\mathcal{D}^yv(\xi^\varepsilon_s)ds.
$$
Applying the It\^o formula to $v(\xi^\varepsilon_t)$, we derive a
new presentation for $U^\varepsilon_t$:
\begin{eqnarray*}
U^\varepsilon_t&=&\frac{1}{\sqrt{\varepsilon}}\int_0^t\nabla^*_z
v(\xi^\varepsilon_s)\sigma(\xi^\varepsilon_s,Y^\varepsilon_s)dB_s
+\frac{1}{\varepsilon}\int_0^t\Big(\mathcal{L}^y-\mathcal{D}^y\Big)
v(\xi^\varepsilon_s)ds
\\
&=&\frac{1}{\sqrt{\varepsilon}}\int_0^t\nabla^*_z
v(\xi^\varepsilon_s)\sigma(\xi^\varepsilon_s,Y^\varepsilon_s)dB_s
-\frac{1}{2\varepsilon}\int_0^t\|\nabla_z v(\xi^\varepsilon_s)
\sigma(\xi^\varepsilon_{s},Y^\varepsilon_s)\|^2ds
\\
&=&\mbox{``continuous martingale$-\frac{1}{2}$ of quadratic
variation process''}.
\end{eqnarray*}
The latter provides that $
Z_t^\varepsilon=\exp\big(U^\varepsilon_t\big) $ is then a positive
continuous local martingale, $Z^\varepsilon_0=1$. Hence, by
Problem 1.4.4. in \cite{LSMar}, $Z_t^\varepsilon$ is a
supermartingale as well, so that $E Z_T^\varepsilon\le 1$.

Set $ \mathfrak{A} =\big\{\int_0^TI\big(\|\xi^\varepsilon_t)\|\ge
n\big)dt>\eta\big\}$. Obviously,
$
1\ge E I_{\mathfrak{A}} Z_T^\varepsilon
$
and this inequality remains valid, if
we replace $Z^\varepsilon_T$ by its lower bound on $\mathfrak{A}$.
Below, we find an appropriate lower bound. Set
$$
\mathcal{A}^y=\sum_{i,j=1}^{d}(\sigma\sigma^*)_{ij}(z,y)\frac{\partial^2}
{\partial z_i\partial z_j}.
$$
Then, we have (see \eqref{Ly}) $ \mathcal{L}^y=
\frac{1}{2}\mathcal{A}^y
+\sum_{i=1}^{d}b_i(z,y)\frac{\partial}{\partial z_i} $ and so, by
\eqref{Dy},
\begin{eqnarray*}
\mathcal{D}^yv(z)=\frac{(2+\|z\|)}{(1+\|z\|)^2}\lef z,b(z,y)\rig+
\frac{1}{2}\mathcal{A}^yv(z) +\frac{1}{2}\|\nabla_z
v(z)\sigma(z,y)\|^2.
\end{eqnarray*}
The function $ \frac{1}{2}\mathcal{A}^yv(z) +\frac{1}{2}\|\nabla_z
v(z)\sigma(z,y)\|^2 $ is bounded, since $\nabla_z v$ and
$\partial^2v$ are bounded and by (A$_a$), $\|\nabla_z
v(z)\sigma(z,y)\|^2$ is bounded too. On the other hand, by (A$_b$)
$$
\lim_{\|z\|\to\infty}\sup_y \frac{\lef
z,b(z,y)\rig}{\|z\|}=-\infty.
$$
Hence, there is a positive constant $K$ such that
\begin{equation}\label{r}
\mathcal{D}^yv(z)\le K\quad\text{and}\quad
\lim_{r\to\infty}\inf_{\|z\|>r}\big(K-\mathcal{D}^yv(z)\big)=\infty.
\end{equation}
We express the lower bound for $Z^\varepsilon_T$ in terms of $K$
and $\mathcal{D}^yv(z)-K$. Taking into account that
$K-\mathcal{D}^yv(z)\ge 0$, write
\begin{eqnarray*}
\log
Z^\varepsilon_T&=&v(\xi^\varepsilon_T)-v(z_0)-\frac{1}{\varepsilon}
\int_0^T\Big(K+\big[\mathcal{D}^yv(\xi^\varepsilon_s)-K\big]\Big)ds
\\
&\ge& -v(z_0)-\frac{K}{\varepsilon}T+
\frac{1}{\varepsilon}\int_0^TI\big(\|\xi^\varepsilon_s\|>n\big)
\big(K-\mathcal{D}^yv(\xi^\varepsilon_s)\big)ds
\\
&\ge& -v(z_0)-\frac{K}{\varepsilon}T+
\frac{1}{\varepsilon}\inf_{\|z\|>r} \big(K-\mathcal{D}^yv(z)\big)
\int_0^TI\big(\|\xi^\varepsilon_s\|>n\big)ds.
\end{eqnarray*}
Set
$$
\log Z_*=-v(z_0)-\frac{K}{\varepsilon}T+
\frac{\eta}{\varepsilon}\inf_{\|z\|>r}
\big(K-\mathcal{D}^yv(z)\big).
$$
Owing to $ \int_0^TI\big(\|\xi^\varepsilon_s\|>n\big)ds>\eta $
on $\mathfrak{A}$, we have $Z^\varepsilon_T\ge Z_*$ on
$\mathfrak{A}$ and, thus, $1\ge E I_{\mathfrak{A}}Z_*$. Now, notice that
$Z_*$ is nonrandom number, so that
$$
\varepsilon^{1-2\kappa}\log P \big(\mathfrak{A}\big)\le
\varepsilon^{1-2\kappa}
v(z_0)+\frac{1}{\varepsilon^{2\kappa}}\Big(K-\eta
\inf_{\|z\|>r}\big(K-\mathcal{D}^yv(z)\Big).
$$
By \eqref{r}, there exists  $r^\circ>0$ such that $
K<\eta\inf_{\|z\|>r^\circ}\big(K-\mathcal{D}^yv(z)). $

Then,
$ \lim_{\varepsilon\to 0}\varepsilon^{1-2\kappa}\log P
\big(\mathfrak{A}\big) =-\infty$.
\end{proof}

\appendix
\section{Stochastic exponential\\ and exponential estimates}
\label{App-B}

The lemma below is borrowed from \cite{LSpok1}. We assume that
$M=(M_t)_{t\ge 0}$ is continuous local martingale with $M_0=0$ and
the predictable variation process $\langle M\rangle_t$. Assume
that $M$ is defined on some stochastic basis
$(\Omega,\mathcal{F},\mathbf{F}=(\mathcal{F}_t)_{t\ge 0},P)$ with
general conditions (see, e.g. \cite{LSMar}, Ch. 1, \S 1). With
$\lambda\in \mathbb{R}$, let us introduce a positive local
martingale
$$
Z_t(\lambda)=\exp\Big(\lambda M_t-\frac{\lambda^2}{2}\langle
M\rangle_t\Big).
$$
It is well known (see e.g. Problem 1.4.4 in \cite{LSMar}) that
$Z_t(\lambda) $ is a supermartingale too. So, owing to
$Z_0(\lambda)=1$, for any stopping time $\tau$,
\begin{equation}\label{B.1}
EZ_\tau(\lambda)\le 1.
\end{equation}

\begin{lemma}\label{lem-B.1}
Let $\tau$ be a stopping time, $\mathfrak{A}\in\mathcal{F}$, and
$\alpha$, $B$ are positive constants.

Then
$$
\begin{aligned}
&{\rm 1)} \ \mathfrak{A}\cap\big\{M_\tau-\frac{1}{2}\langle
M\rangle_\tau\ge\alpha\big\}= \big\{M_\tau-\frac{1}{2}\langle
M\rangle_\tau\ge\alpha\big\} \Rightarrow P(\mathfrak{A})\le
e^{-\alpha}
\\
&{\rm 2)} \ \mathfrak{A}\cap\big\{M_\tau\ge\alpha, \ \langle
M\rangle_\tau \le B\big\}=\big\{M_\tau\ge\alpha, \ \langle
M\rangle_\tau \le B\big\} \Rightarrow P(\mathfrak{A})\le
e^{-\frac{\alpha^2}{2B}}
\\
&{\rm 3)} \ P(\sup_{t\le T}|M_t|\ge \alpha, \langle M\rangle_T\le
B)\le 2e^{-\frac{\alpha^2}{2B}}.
\end{aligned}
$$
\end{lemma}

\begin{proof}
\mbox{}

1) By \eqref{B.1}, $ 1\ge EI_{\mathfrak{A}}Z_\tau(1)\ge
P(\mathfrak{A})e^\alpha $ and the result holds.

2) By \eqref{B.1}, $1\ge EI_{\mathfrak{A}}Z\tau(\lambda)\ge
P(\mathfrak{A}) e^{\lambda\alpha-\frac{\lambda^2B}{2}}\ge
P(\mathfrak{A})\ge e^{\frac{\alpha^2}{2B}}$ (the latter inequality
under $\lambda=\frac{\alpha}{B}$). Then, the assertion follows.

3) Introduce Markov times $\tau_\pm=\inf\{t:\pm M_t\ge \alpha\}$
(here $\inf\{\emptyset\}=\infty$) and two sets
$\mathfrak{A}_\pm=\{\tau_\pm\le T,\langle M\rangle_T\le B\}$. By
2), $P(\mathfrak{A}_\pm)\le e^{-\frac{\alpha^2}{2B}}$.

The assertion is valid, since  $\{\sup_{t\le T}|M_t|\ge
\eta\}\subseteq\mathfrak{A}_+\cup\mathfrak{A}_-$.
\end{proof}

\section{Exponential negligibility of
$\varepsilon^l\sup_{t\le T}\|\xi_t^\varepsilon\|^p$} \label{App-A}

The next Lemma plays a crucial role in many proofs of this paper and may be
of independent interest.

\begin{lemma}\label{theo-A.1}
Assume $p>0$, and $l>\frac{p}{2}$. Then for any $\eta>0$ and
$T>0$
$$
\lim_{\varepsilon\to 0}\varepsilon^{1-2\kappa}\log
P\Big(\sup_{t\le T}
\varepsilon^{l}\|\xi^\varepsilon_t\|^p>\eta\Big)=-\infty.
$$
\end{lemma}
\begin{proof}
Show first that for any $L>0$
\begin{equation}\label{B.II}
\lim_{\varepsilon\to 0} \varepsilon^{1-2\kappa} \log P\Big(
\varepsilon^{l}\sup_{t\le T}\|\xi_t^\varepsilon\|^p > \eta,
\int_0^T \|\xi_s^\varepsilon\|^2ds \le L \Big)=-\infty.
\end{equation}
Noticing that
$\|\xi^\varepsilon_t\|^2=(\xi^\varepsilon_t)^*\xi^\varepsilon_t$,
by the It\^o formula we find

\begin{equation}\label{hado}
\begin{aligned}
\|\xi_t^\varepsilon\|^2&=
\|z_0\|^2+\varepsilon^{-1}\int_0^t2\lef\xi_s^\varepsilon,
b(\xi_s^\varepsilon,Y^\varepsilon_s)\rig ds
\\
& \quad +\varepsilon^{-1}\int_0^t\trace\big\{
\sigma^*(\xi_s^\varepsilon,Y^\varepsilon_s)\sigma(\xi_s^\varepsilon,Y^\varepsilon_s)
\big\}ds
\\
& \quad +\varepsilon^{-1/2} M^\varepsilon_t,
\end{aligned}
\end{equation}
where
$
M^\varepsilon_t= \int_0^t2\lef
\xi^\varepsilon_s,\sigma(\xi^\varepsilon_s,Y^\varepsilon_s)
dB_s\rig
$
the continuous martingale with
the predictable variation process
\begin{equation}\label{Me}
\langle M^\varepsilon\rangle_t=4\int_0^t(\xi^\varepsilon_s)^*
a(\xi^\varepsilon_s,Y^\varepsilon_s)\xi^\varepsilon_sds.
\end{equation}
By (A$_a$),
$
\trace\big\{
\sigma^*(\xi_s^\varepsilon,Y^\varepsilon_s)\sigma(\xi_s^\varepsilon,Y^\varepsilon_s)
\big\}
$
is bounded. By (A$_b$), $\lef
\xi_s^\varepsilon,b(\xi_s^\varepsilon,Y^\varepsilon_s)\rig$ is
negative for $\|\xi_t^\varepsilon\|>C$ and $C$ large enough, so
that $
\lef\xi_s^\varepsilon,b(\xi_s^\varepsilon,Y^\varepsilon_s)\rig $
is bounded from above. Hence, with some positive constant $K$,
\eqref{hado} provides
\begin{equation}\label{B.M}
\sup_{t\le T}\|\xi_t^\varepsilon\|^2\le
\varepsilon^{-1}K+\varepsilon^{-1/2}\sup_{t\le T}|M_t^\varepsilon|.
\end{equation}
Further, since $\sigma^*\sigma$ is nonnegative definite matrix, by (A$_a$)
\begin{equation}\label{A.5m}
\langle M^\varepsilon\rangle_T \le 4\Lambda
\int_0^T\|\xi^\varepsilon_s\|^2ds.
\end{equation}
We use now \eqref{B.M} and \eqref{A.5m} for the proof of \eqref{B.II}.
From \eqref{B.M}, with the help of H\"older's inequality for $p>1$ and simply $(a+b)^p\le a^p+b^p$ for $p\le 1$ and $a,b>0$, we derive (recall that
$2l>p$)
\begin{equation}\label{A.4nn}
\varepsilon^{2l}\sup_{t\le T}\|\xi_t^\varepsilon\|^{2p}\le (2^{p-1}\vee 1)\left(
K^p\varepsilon^{2l-p}+\varepsilon^{2l-p/2}\sup_{t\le T}|M_t^\varepsilon|^p\right).
\end{equation}
Then, by \eqref{A.4nn} and \eqref{A.5m}, we have
$$
\begin{aligned}
&\Big\{\varepsilon^{2l}\sup_{t\le
T}\|\xi_t^\varepsilon\|^{2p}>\eta^2,
\int_0^T\|\xi_s^\varepsilon\|^2ds\le L\Big\}
\\
&\quad\subseteq \left\{\sup_{t\le T}|M_t^\varepsilon|\ge
\Bigg(\frac{\eta^2/(2^{p-1}\vee 1)-K^p\varepsilon^{2l-p}}{\varepsilon^{2l-p/2}}\Bigg)^{1/p},
\langle
M^\varepsilon\rangle_T \le \frac{L}{4\Lambda}\right\}.
\end{aligned}
$$
With $\varepsilon_0$ such that
$
\eta^2/(2^{p-1}\vee 1)-K^p\varepsilon^{2l-p}_0=\eta^2/(2^p\vee 2),
$
and  $\varepsilon<\varepsilon_0$ the above inclusion
provides the following inequality:
$$
\begin{aligned}
& P\Big(\varepsilon^{l}\sup_{t\le T}\|\xi_t^\varepsilon\|^p>\eta,
\int_0^T\|\xi_s^\varepsilon\|^2ds\le L\Big)
\\
&\quad \le P\Big(\sup_{t\le T}|M_t^\varepsilon|\ge
\frac{\eta^{2/p}}{(2\vee 2^{1/p})\varepsilon^{(2l-p/2)(1/p)}},
\langle M^\varepsilon\rangle_T \le \frac{L}{4\Lambda}\Big).
\end{aligned}
$$
Now, by Lemma \ref{lem-B.1}(3) (Appendix \ref{App-B}), we find
that
$$
P\Bigg(\varepsilon^{l}\sup_{t\le T}\|\xi_t^\varepsilon\|^p>\eta,
\int_0^T\|\xi_s^\varepsilon\|^2ds\le L\Bigg)\le 2\exp\left(-
\frac{4\Lambda\eta^{4/p}}{2(4\vee 2^{2/p})L\varepsilon^{(2l-p/2)(2/p)}}.
\right)
$$
Since $l>\frac{p}{2}$, it holds $(2l-p/2)(2/p)\ge 1$. Hence,
\begin{gather*}
\varepsilon^{1-2\kappa}\log P\Big(\varepsilon^{l}\sup_{t\le
T}\|\xi_t^\varepsilon\|^p>\eta,
\int_0^T\|\xi_s^\varepsilon\|^2ds\le L\Big)
\\
\le \text{const.}\left(\varepsilon^{1-2\kappa}-\frac{1}{\varepsilon^{2\kappa}}\right)
\to-\infty, \ \varepsilon\to 0.
\end{gather*}

In the next step of proof, we show that for $L$ large enough
\begin{equation}\label{A.3a}
\lim_{\varepsilon\to 0}\varepsilon^{1-2\kappa}\log
P\left(\int_0^T\|\xi_s^\varepsilon\|^2ds> L\right)=-\infty.
\end{equation}
With $\inf\{\emptyset\}=\infty$, set $ \tau^\varepsilon=\inf\{t\ge
0:\int_0^t\|\xi^\varepsilon_s\|^2ds\ge L\} $ and notice that
\eqref{A.3a} holds, if for $L$ large enough
\begin{equation}\label{A.5a}
\lim_{\varepsilon\to 0}\varepsilon^{1-2\kappa}\log P
\big(\tau^\varepsilon\le T\big)=-\infty.
\end{equation}
To this end, let us evaluate evaluate from below $M_{\tau^\varepsilon}:=
\varepsilon^{1/2}M^\varepsilon_{\tau^\varepsilon}$.
By \eqref{hado}, we have
$$
\begin{aligned}
M_t&\ge -\varepsilon\|z_0\|^2-
\int_0^t2\big|\lef \xi_s^\varepsilon,
b(\xi_s^\varepsilon,Y^\varepsilon_s)\rig\big|
I(\|\xi^\varepsilon_s\|\le C) ds
\\
&\quad -\int_0^t2\lef \xi_s^\varepsilon,
b(\xi_s^\varepsilon,Y^\varepsilon_s)\rig
I(\|\xi^\varepsilon_s\|>C) ds
\\
& \quad -\int_0^t\trace\big\{
\sigma^*(\xi_s^\varepsilon,Y^\varepsilon_s)\sigma(\xi_s^\varepsilon,Y^\varepsilon_s)
\big\}ds.
\end{aligned}
$$
From (A$_a$) and (A$_b$), it follows the existence of a positive constants
$K$ such that

$$
\begin{aligned}
&-\varepsilon\|z_0\|^2-
\int_0^T2\big|\lef \xi_s^\varepsilon,
b(\xi_s^\varepsilon,Y^\varepsilon_s)\rig\big|
I(\|\xi^\varepsilon_s\|\le C) ds
\\
&\quad -\int_0^T2\lef \xi_s^\varepsilon,
b(\xi_s^\varepsilon,Y^\varepsilon_s)\rig
I(\|\xi^\varepsilon_s\|>C) ds
\\
& \quad -\int_0^T\trace\big\{
\sigma^*(\xi_s^\varepsilon,Y^\varepsilon_s)\sigma(\xi_s^\varepsilon,Y^\varepsilon_s)
\big\}ds
\\
&\ge -K+r\int_0^{\tau^\varepsilon}\|\xi^\varepsilon_s\|^2I(\|\xi^\varepsilon_s\|>C)ds
\\
&\ge -(K+rC^2T)+r\int_0^{\tau^\varepsilon}\|\xi^\varepsilon_s\|^2ds.
\end{aligned}
$$
So, on $\{\tau^\varepsilon\le T\}$, we have
$
M_{\tau^\varepsilon}\ge -(K+CT)+rL.
$
Since $\langle M\rangle_t=\varepsilon\langle M^\varepsilon\rangle_t$,
by \eqref{Me} and (A$_a$), we have
$
\langle M\rangle_{\tau^\varepsilon}\le
4\Lambda\varepsilon\int_0^{\tau^\varepsilon}\|\xi^\varepsilon_s\|^2ds=
4\Lambda L\varepsilon.
$

Consequently, by Lemma \ref{lem-B.1}(3), with $L$ so large that
$-(K+CT)+rL=1$, we find that
$
P\left(\tau^\varepsilon\le T\right)\le 2e^{-1/(\varepsilon 8\Lambda L)},
$
that is
$$
\varepsilon^{1-2\kappa}\log P\left(\tau^\varepsilon\le T\right)\le
\varepsilon^{1-2\kappa}\log 2-\frac{1}{8\Lambda L\varepsilon^{2\kappa}}\to-\infty, \
\varepsilon\to 0.
$$
\end{proof}

\section{Exponential negligibility\\ for semimartingales}
\label{sec-7}

Let $S^\varepsilon_t$ and $N^\varepsilon_t$ be continuous
semimartingale and martingale respectively with paths in
$\mathbb{C}_{[0,T]}(\mathbb{R}^\ell)$. Set
$$
Y^\varepsilon_t=S^\varepsilon_t+\varepsilon^{1/2-\kappa}N^\varepsilon_t.
$$
Denote $S^\varepsilon_t(i)$, $N^\varepsilon_t(i)$ entries of
$S^\varepsilon_t$, $N^\varepsilon_t$  and
$Y^\varepsilon_t(i)=S^\varepsilon_t(i)+\varepsilon^{1/2-\kappa}N^\varepsilon_t(i)$.
Define $|S^\varepsilon_t|=\sum_{i=1}^\ell|S^\varepsilon_t(i)|$ and
similarly $|N^\varepsilon_t|$, $|Y^\varepsilon_t|$.

\begin{lemma}\label{lem-A}
Assume for some nonnegative $c_1\div c_3$
\[
\begin{aligned}
& |S^\varepsilon_t|\le
c_1+c_2\int_0^t\big(1+|Y^\varepsilon_s|\big)ds
\\
& d\langle N^\varepsilon(i)\rangle_t\le c_3dt, \ i=1\ldots,\ell.
\end{aligned}
\]
Then for any $T>0$
$$
\lim_{C\to\infty}\varlimsup_{\varepsilon\to
0}\varepsilon^{1-2\kappa}\log P\Big( \sup_{t\le
T}|Y^\varepsilon_t|>C\Big)=-\infty.
$$
\end{lemma}

The proof can be found in \cite{LSV} (Lemma A.1).

\end{document}